\documentclass[11pt,reqno]{amsart}
\usepackage{amsmath}
\usepackage{amsxtra}
\usepackage{amscd}
\usepackage{amsthm}
\usepackage{amsfonts}
\usepackage{amssymb}
\usepackage{eucal}

\newtheorem{theorem}{Theorem}[section]

\newtheorem{cor}[theorem]{Corollary}
\newtheorem{lem}[theorem]{Lemma}
\newtheorem{prop}[theorem]{Proposition}

\newcommand{\bean}{\begin{eqnarray}}
\newcommand{\eean}{\end{eqnarray}}
\newcommand{\be}{\begin{displaymath}}
\newcommand{\ee}{\end{displaymath}}
\newcommand{\bea}{\begin{eqnarray*}}   
\newcommand{\eea}{\end{eqnarray*}}

\newcommand{\Ref}[1]{{$($\ref{#1}$)$}}

\newcommand{\nc}{\newcommand}
\nc{\on}{\operatorname}
\nc{\ch}{\mbox{ch}}
\nc{\Z}{{\mathbb Z}}
\nc{\C}{{\mathbb C}}
\nc{\pone}{{\mathbb C}{\mathbb P}^1}
\nc{\pa}{\partial}
\nc{\al}{\alpha}
\nc{\ri}{\rangle}
\nc{\lef}{\langle}
\nc{\la}{\lambda}
\nc{\ep}{\epsilon}

\nc{\mc}{\mathcal}

\begin{document}

\title[Diagonal shifted by roots of unity]
{Symmetric polynomials vanishing on the diagonals shifted by roots of unity.}

\author{B. Feigin, M. Jimbo, T. Miwa, E. Mukhin and Y. Takeyama}
\address{BF: Landau institute for Theoretical Physics, Chernogolovka,
142432, Russia}\email{feigin@feigin.mccme.ru}  
\address{MJ: Graduate School of Mathematical Sciences, University of
Tokyo, Tokyo 153-8914, Japan}\email{jimbomic@ms.u-tokyo.ac.jp}
\address{TM: Department of Mathematics, Graduate School of Science, 
Kyoto University, Kyoto 606-8502 Japan}\email{tetsuji@kusm.kyoto-u.ac.jp}
\address{EM: Department of Mathematics, 
Indiana University-Purdue University-Indianapolis, 
402 N.Blackford St., LD 270, 
Indianapolis, IN 46202}
\email{mukhin@math.iupui.edu}
\address{YT: Department of Mathematics, Graduate School of Science, 
Kyoto University, Kyoto 606-8502 Japan}
\email{takeyama@kusm.kyoto-u.ac.jp} 

\date{\today}

\setcounter{footnote}{0}\renewcommand{\thefootnote}{\arabic{footnote}}

\begin{abstract}
For a pair of positive integers $(k,r)$ with $r\geq 2$ such that $k+1$
and $r-1$ are  
relatively prime, we describe the space of symmetric polynomials in variables
$x_1,\dots,x_n$ which vanish at all 
diagonals of codimension $k$ of the form $x_i=tq^{s_i}x_{i-1}$, 
$i=2,\dots,k+1$, where $t$ and $q$ are primitive roots of unity of
orders $k+1$ and $r-1$.
\end{abstract}
\maketitle

\renewcommand{\thefootnote}{\arabic{footnote})}
\renewcommand{\arraystretch}{1.2}


\section{Introduction: the wheel condition}
Fix a finite set $\mc S$ of non-zero complex numbers called the
{\it wheel set}. 
A symmetric polynomial $f\in \C[x_1,\cdots,x_n]^{{\mathfrak S}_n}$ 
satisfies the {\it wheel condition} 
relative to $\mc S$ if
$f(x_1,\dots,x_n)=0$ on all planes 
which have the form
\bean\label{wheel}
x_2=t_1x_1,\;\; x_3=t_2x_2,\;\;\dots,\;\; 
x_l=t_{l-1}x_{l-1},\;\; x_1=t_lx_l,
\eean
where $t_1,\cdots,t_l\in\mc S$. 
Note that \Ref{wheel} implies that $\prod_{i=1}^lt_i=1$.
This condition is called the resonance condition.
We denote the space of all symmetric polynomials satisfying the wheel
condition associated to $\mc S$ 
by $F_{\mc S}$. 

Note that the space of symmetric polynomials in $n$ variables
satisfying the wheel condition 
relative to ${\mc S}$ is an ideal in $\C[x_1,\dots,x_n]$.

Let $g_{nm}$ be the dimension of the symmetric polynomials in $n$
variables of degree $m$ satisfying the wheel condition. Then the 
character is given by
\be
\chi(F_{\mc S}^k)=\sum_{n,m} g_{nm}z^nv^m.
\ee
The basic question we are interested in is the computation of the
character and a construction of an explicit basis in $F_{\mc S}$.

The study of the polynomials satisfying the wheel condition was
initiated in \cite{FJMM2}. We briefly describe the results of that paper.

For natural numbers $k,r$, $r\geq 2$, we fix $t,q\in \C$ such that
the resonance condition $q^at^b=1$ is valid
if and only if $a=(k+1)s$, $b=(r-1)s$ for some $s$.
Define the wheel set $\mc S_r(q,t)\subset \C$ by 
\bean\label{set}
\mc S_r(q,t)=\{t,tq,\dots,tq^{r-1}\}.
\eean
The number of variable $l$ in (\ref{wheel}) in this case, can be
a multiple of $k+1$. However, it is easy to see that it is enough to impose
the zero condition only for $l=k+1$. This remark is also valid in the root of
unity case which we will consider in this paper.

It is proved in \cite{FJMM2} that the space $F_{\mc S_r(q,t)}$ has a
basis of Macdonald polynomials $\{P_\la(x_1,\dots,x_n;q,t)\}$ 
where $\la$ ranges over the $(k,r,n)$-admissible partitions, satisfying
\be
\la_i-\la_{i+k} \geq r \qquad (i=1,\dots,n-k).
\ee
In particular, the Macdonald polynomials corresponding to the
admissible partitions are well defined. Admissible partitions appeared
first in the work \cite{P}, a bosonic formula for the
corresponding character is given in \cite{FJLMM}.

The case $k+1$ and $r-1$ are relatively prime is of special
interest. In this case we can assume
\be
t=u^{r-1}, \qquad q=u^{-(k+1)},
\ee
for some $u$ which is not a root of unity. 

The Jack limit $u\to 1$ of the 
space $F_{\mc S_r(q,t)}^k$ is spanned by the set of Jack polynomials
 $\{J_\la(x_1,\dots,x_n;\beta)\}$ where $\beta=-(r-1)/(k+1)$ and $\la$
 ranges over $(k,r,n)$-admissible partitions. If $(k+1,r-1)=1$ then it
 is expected that the limiting space coincides
with the space of correlation
  functions of an abelian current in a vertex operator algebra,
  associated with the minimal series $(k+1,k+r)$ of the $W_k$ algebra,
  see \cite{FJLMM}. 

>From now on we assume $k+1$ and $r-1$ are relatively prime.
The present paper could be viewed as the ``root of unity'' version of
\cite{FJMM2}. Namely, we consider the case when $u$ is a primitive
root of unity of order $(r-1)(k+1)$. We have then
\bean\label{roots}
t^{k+1}=1\quad(\hbox{a primitive root}),\qquad q^{r-1}=1\quad(\hbox{a primitive root}).
\eean
We denote the space of symmetric function satisfying the wheel
condition related to 
the wheel set $S_{r-1}(q,t)$ with \Ref{roots}
by $F^{(k,r)}$ and the cooresponding character by $\chi_{k,r}(z,v)$. 
We also denote by $F^{(k,l)}_n$ the subspace of $F^{(k,l)}$ consisting of
the functions of $n$ variables. Note that $q^{r-1}=1$, and the new wheel
set with $r-1$ is thus related to the previous one $S_r(q,t)$.
However, in the present case, it is not true that the resonance condition
$q^at^b=1$ necessarily restricts to the case
$a=(k+1)s$, $b=(r-1)s$ for some $s$.
In other words, not all vanishing planes \Ref{wheel} at the roots of unity
case are limits of vanishing planes for generic $u$. Therefore,
we have new phenomenon in the root of unity case, which is independent
of generic results of \cite{FJMM2}.

The plan of paper is as follows. In Section 2 we prepare
some special polynomials from MacDonald's book \cite{M}.
In Section 3 we give a basis of the space $F^{(k,2)}$.
In Scetion 4 we give the main results of this paper,
Theorems \ref{basis thm} and \ref{char thm},
which give an explicit basis in $F^{(k,r)}$ and
the character $\chi_{k,r}(z,v)$, respectively.



\section{Macdonald polynomials} In this section we recall some facts
we need about the theory of Macdonald polynomials and fix our notations.

\subsection{Partitions}
We denote $\pi_n$ the set of partitions of length at most $n$.
By this we mean that an element $\lambda$ of $\pi_n$ is
a nonincreasing sequence of nonnegative integers
$\la=(\la_1,\dots,\la_n)$, $\la_1\geq \dots \geq \la_n\geq 0$.
The sum $|\la|=\sum_i\la_i$ is called weight of $\la$.
We denote the number of parts $\la_j$ which are equal to $i$ by
$m_i=m_i(\la)$ and write $\la=(0^{m_0},1^{m_1},2^{m_2},\dots)$. 

 %
A partition $\la\in\pi_n$ is called {\it $(k,r,n)$-admissible} if
$\la_i-\la_{i+k}\geq r$ for $i=1,\dots,n-k$. We will use
$(k,1,n)$-admissibility in this paper. This is equivalent to the condition
that $m_i(\lambda)\leq k$ for all $i$. In the below we simply call
$\lambda$ is (non)-admissible with the understanding that it means the
(non)-admissibility for $(k,1,n)$ with a prescribed value of $k$.

There is a partial {\it dominance} order of
partitions. For two 
partitions $\la,\mu$ of the same weight we write
$\la>\mu$ iff $\la\neq \mu$, $\la_1+\dots+\la_i\geq \mu_1+\dots+\mu_i$
for all $i$. 

We also use the total {\it lexicographical} order of partitions. For two 
partitions $\la,\mu$ of the same weight we write $\la\prec\mu$ iff
for some $i$ we have $\la_i<\mu_i$ and $\la_j=\mu_j$, $j=1,\dots,i-1$.

\subsection{Macdonald polynomials}
The Macdonald operators $D_n^{r}(q,t),$ $0\le r\le n$, 
are mutually commuting $q$-difference operators acting on 
the ring of symmetric polynomials $\C(q,t)[x_1,\cdots,x_n]^{{\mathfrak S}_n}$:
\be\label{Dr}
D_n^r=\sum_{|I|=r}A_I(x;t)T_{I},
\ee
where 
\bea
&&A_I(x;t)=t^{r(r-1)/2}\prod_{i\in I \atop j\not\in I}
\frac{tx_i-x_j}{x_i-x_j}, \\
&&T_{I}=\prod_{i\in I}T_{q,x_i},
\quad 
(T_{q,x_i}f)(x_1,\cdots,x_n)=f(x_1,\cdots,qx_i,\cdots,x_n),
\eea
and $I\subset\{1,\cdots,n\}$ runs over subsets of cardinality $r$. 
Let $D_n(X;q,t)=\sum_{r=0}^{n}D_n^rX^r$ be their generating function.

For a partition $\la\in\pi_n$ the 
{\it Macdonald polynomial} $\{P_\lambda\}$ is defined as
a unique eigenvector of $D_n(X;q,t)$ of the form
\be
P_\lambda=m_\lambda+\sum_{\mu<\lambda}u_{\lambda\mu}m_\mu 
\qquad (u_{\lambda\mu}\in \C(q,t)),  
\ee
where
\be
m_\la= \prod_{i=1}^n\frac1{m_i!}\sum_{w\in{\mathfrak S}_n}w\left(
  x_1^{\la_1}\dots x_n^{\la_n}\right)
\ee
is the monomial symmetric function. Here a permutation $w\in{\mathfrak S}_n$
acts on a function by permuting the arguments
$wf(x_1,\dots,x_n)=f(x_{w^{-1}(1)},\dots, x_{w^{-1}(n)})$. 

Such a polynomial is unique and the coefficients $u_{\lambda\mu}$ 
are rational functions of $q$ and $t$ with possible poles of the form
\be\label{pole}
q^at^b=1\qquad (a\in\Z_{\ge 0},b\in\Z_{>0}),  
\ee
see \cite{Int}.

The corresponding eigenvalues are given by the formula:
\be\label{diag}
D_n(X;q,t)P_\lambda=\prod_{i=1}^n(1+Xq^{\lambda_i}t^{n-i}) P_\lambda,
\ee
see \cite{M}, VI,3,4.

Note that the degree of $P_\la$ is equal to the weight of $\la$.

\subsection{Hall-Littlewood polynomials} 

Let $\la$ be a partition of length $n$ and let $t\in \C$. The {\it
  Hall-Littlewood polynomial} $P_\la(x_1,\dots,x_n;t)$ is defined by
\bean\label{HL}
P_\la(x_1,\dots,x_n;t)
=\prod_{i=1}^n\prod_{j=1}^{m_i}\frac{1-t}{1-t^j}\sum_{w\in{\mathfrak S}_n}
w\left( x_1^{\la_1}\dots x_n^{\la_n} \prod_{1\leq i<j\leq n}
  \frac{x_i-tx_j}{x_i-x_j}\right).
\eean
We will construct a basis of $F^{(k,2)}_n$ by using these polynomials
in Proposition \ref{r=2 prop}.

It is well known that for generic $t$ the function $P_\la$ is a symmetric
polynomial of the form
\bean\label{tri}
P_\la=m_\la+\sum_{\mu<\la} a_{\la\mu}m_\mu, \qquad a_{\la\mu}\in \C,
\eean
 see \cite{M}, III 1,2. 

It is also well known that the Hall-Littlewood polynomials are
obtained from the Macdonald polynomials by setting $q=0$:
$P_\la(x_1,\dots,x_n;t)=P_\la(x_1,\dots,x_n;0,t)$, see \cite{M}, p.324.
Although for $r=2$ we are dealing with the case $q=1$ in the wheel
condition, the solution is given in terms of the special case of MacDonald
polynomials with $q=0$. This is not surprising because as we have already
remarked the root of unity case is not a specialization of the generic case.
(See Proposition \ref{r=2 general}.)

\section{The case $r=2$} 

\subsection{Preliminaries}
Let $r=2$. 
Then we have $q=1$, $t$ a primitive root of unity of order $k+1$,
the wheel set is $\{t\}$, and $F^{(k,2)}$ is the space of symmetric
functions satisfying a single condition:
\bean\label{r=2}
f(x,tx,t^2x,\dots,t^kx,x_{k+2},\dots,x_n)=0.
\eean

The condition \Ref{r=2} appears as the initial condition for the
recurrence relation of ``deformed cycles'' in the 
construction of form factors in
$SU(k+1)$ invariant Thirring model with values in tensor products of
vector representations, see Proposition 7.2 in \cite{T}. 

Consider the example $k=1$ and $r=2$. We have $t=-1$, and
the space $F^{(1,2)}_n$ consists of polynomials of the form
\be
\prod_{1\leq i<j \leq n}(x_i+x_j)g(x_1,\dots,x_n),
\ee
where $g$ is any symmetric function. In the generic case
for $k=1$ and $r=2$ the wheel set is $\{t,t^{-1}\}$ with $t\not=-1$.
The functions satisfying the wheel condition
have the form
\be
\prod_{1\leq i< j \leq n}(x_i-tx_j)(x_j-tx_i)g(x_1,\dots,x_n).
\ee
The space $F^{(1,2)}_n$ is
greater than the space of functions for the generic case in the limit $t=-1$.

\subsection{The space $E^{(k,2)}$}
We describe the dual space $E^{(k,2)}_n=(F_n^{(k,2)})^*$ to $F_n^{(k,2)}$. 

Let $F_n=\C[x_i]^{{\mathfrak S}_n}_{i=1,\dots,n}$
be a commutative algebra of symmetric 
polynomials in variables $x_1,\ldots,x_n$. We have $F_n^{(k,2)}\subset F_n$. 
The space $F_n$ has a basis of elementary monomial functions 
$\{m_\la\}_{\la\in\pi_n}$. We
use the usual degree, i.e., ${\rm deg}x_i=1$ and the degree
of $m_\la$ is equal to the weight of partition $\la$.

Let $E=\C[e_i]_{i\in\Z_{\geq0}}$ be a commutative algebra of polynomials
in variables $e_i$. 
For a partition $\la=(\la_1,\dots,\la_n)\in\pi_n$,  
$\la_1\geq \dots\geq \la_n$, set $e_\la=\prod_{i=1}^ne_{\la_i}$. 
We define the degree of $e_i$ to be $i$. Then, the degree of
$e_\la$ is equal to the weight of $\la$.
For $\mu\in\pi_m$ and $\nu\in\pi_n$ we define $\lambda=\mu\cup\nu\in\pi_{m+n}$
by joining the parts of $\mu$ and $\nu$. Then, we have $e_\lambda=e_\mu e_\nu$.

Monomials $e_\la\in\pi_n$
are linearly independent. Let $E_n\subset E$ be the subspace with the 
basis $\{e_\la\}_{\la\in\pi_n}$. 

Let 
\be
e(z)=\sum_{i\geq 0}e_iz^i,  
\ee
be a formal power series in $z$. Then $E_n$ is spanned by the 
coefficients of power series $\prod_{i=1}^n e(z_i)$. Moreover, we have
\be
e(z_1)\dots e(z_n)=\sum_{\la\in\pi_n} e_\la m_\la(z_1,\dots,z_n).
\ee 

Define a bilinear pairing
$\langle, \rangle:\; E_n\otimes F_n\to \C$ by setting
\be
\langle
e(z_1)\dots e(z_n),f(x_1,\dots,x_n)\rangle=f(z_1,\dots,z_n),
\ee
for any symmetric polynomial $f(x_1,\dots,x_n)\in F_n$.

The following lemma is clear.
\begin{lem}
The pairing $\langle,\rangle$ is a well defined bilinear nondegenerate
pairing of graded spaces. Moreover, the bases $e_\la$ and $m_\la$ are dual, 
$\langle e_\la,m_\mu\rangle=\delta_{\la\mu}$. $\square$
\end{lem}

Let $t$ be a primitive root of unity of order $k+1$.
We have a subspace $F_n^{(k,2)}\subset F_n$. 
Let $J$ be the space spanned by the coefficients of formal power series
$\prod_{i=0}^ke(t^iz)$.

\begin{lem}
The orthogonal complement $(F_n^{(k,2)})^\perp\subset E_n$ with respect to the
pairing $\langle,\rangle$ coincides with the subspace $J\cdot E_{n-k-1}$,
\be 
(F_n^{(k,2)})^\perp=J\cdot E_{n-k-1}.
\ee
\end{lem}
\begin{proof} 
For $f\in F_n$, we have
\be
\langle \prod_{i=0}^ke(t^iz)\prod_{i=k+2}^n
e(z_i),f(x_1,\dots,x_n)\rangle= f(z,tz,\dots,t^kz,z_{k+2},\dots,z_n)=0,
\ee
if and only if $f\in F_n^{(k,2)}$. Therefore $(J\cdot
E_{n-k-1})^\perp=F_n^{(k,2)}$. Since the graded components are 
finite-dimensional and the pairing respects the grading, we obtain $(J\cdot
E_{n-k-1})^{\perp\perp}=J\cdot E_{n-k-1}$ and the lemma is proved.
\end{proof}

Denote by $E_n^{(k,2)}$ the quotient space of $E_n$ by the
space $J\cdot E_{n-k-1}$,
\be
E_n^{(k,2)}=E_n/\left(J\cdot E_{n-k-1}\right).
\ee
Since the space of relations is graded, the 
subspace $E_n^{(k,2)}$ inherits grading from $F_n$.

The following is clear.

\begin{lem}\label{dual lem} 
The pairing $\langle,\rangle$ induces a well defined nondegenerate 
bilinear pairing of graded spaces
\be
\langle,\rangle:\; E_n^{(k,2)}\otimes F_n^{(k,2)} \to \C.
\ee
\qquad\qquad \qquad\qquad\qquad\qquad\qquad\qquad
\qquad\qquad\qquad\qquad\qquad\qquad\qquad\qquad
$\square$
\end{lem}

\subsection{The spanning set of $E_n^{(k,2)}$} 
Denote the series $\prod_{i=0}^ke(t^iz)$ by $\mc E_k(z)$. We have 
$\mc E_k(tz)=\mc E_k(z)$ and therefore $\mc E_k(z)$ has the form
\be
\mc E_k(z)=\sum_{i\geq0} \ep_iz^{i(k+1)} \qquad (\ep_i\in E_{k+1}). 
\ee

We have
\bean\label{ep}
\ep_i=(-1)^{ik}e_i^{k+1}+\sum_{\la\succ(i^{k+1})} c_{\la,i}e_\la
\quad(c_{\la,i}\in \C). 
\eean 
Note that $e_i^{k+1}=e_{(i^{k+1})}$.

We denote $\bar e$ the image of an element $e\in E_n$ 
in the quotient space $E_n^{(k,2)}$.

\begin{lem}\label{span lem}
The elements $\{\bar e_\la\}$ with admissible partitions $\la\in\pi_n$,
i.e., those satisfying the condition $m_i(\la)\leq k$ for all $i$,
span $E_n^{(k,2)}$. 
\end{lem} 
\begin{proof}
We fix nonnegative integers $d,n$ and
work with a finite-dimensional space generated by $e_\la$ with
partitions $\la\in\pi_n$ of weight $d$. 

We claim that if $\la$ is nonadmissible, i.e., $m_i>k$ for some $i\geq 0$,
then $\bar e_\la$ is a linear combination of $\bar e_\nu$
such that $\nu\succ\la$.

Indeed, if $i=0$ then $\bar e_\la=0$ and
there is nothing to prove. Otherwise let $\tilde \la\in\pi_{n-k-1}$ 
be a partition
obtained from $\la$ by deleting $k+1$ parts equal to $i$. Then we have
$e_\la=e_{(i^{k+1})}e_{\tilde\la}$ where $e_{\tilde\la}\in E_{n-k-1}$. 

We use the relation $\bar \ep_i \bar e_{\tilde\la}=0$ and get (see \Ref{ep})
\be
\bar e_\la=\bar e_{(i^{k+1})}\bar e_{\tilde\la}
=-(-1)^{ik}\sum_{\mu\succ(i^{k+1})} c_{\mu,i}\bar e_{\mu\cup\tilde\lambda}.
\ee
and our claim is proved.

Now the lemma follows. Indeed, we rewrite $\bar e_\la$ in terms of $\bar e_\nu$
with larger $\nu$ (with respect to $\succ$ ordering). Then we
rewrite nonadmissible $\bar e_\nu$ appearing in this sum
in terms of $\bar e_\rho$ with even larger partitions $\rho$ and so on.
Since our space is finite-dimensional, after finitely many
repetitions we obtain a sum with only admissible partitions.
\end{proof}

In fact we will show (see Corollary \ref{basis cor})
that the elements $\{\bar e_\la\}$ with admissible
partitions $\la\in\pi_n$ 
are linearly independent and therefore form a basis in $E_n^{(k,2)}$.

\subsection{A basis in $F^{(k,2)}$} \label{HL basis sec}

Let $t$ be a primitive root of unity of order $k+1$.
\begin{prop}\label{r=2 prop}
The set of Hall-Littlewood polynomials $\{P_\la(x_1,\dots,x_n;t)\}$,
where $\la$ ranges over all admissible partitions is a basis
of $F^{(k,2)}$.
\end{prop}
\begin{proof}
The admissibility of $\lambda$ is nothing but $m_i(\lambda)\leq k$.
It implies that polynomials 
$\{P_\la(x_1,\dots,x_n;t)\}$ have no pole when the variable $t$
in (\ref{HL}) specializes to the root of unity, $t^{k+1}=1$. Because of
\Ref{tri} they are linearly independent. It is also clear from the definition 
that they satisfy the condition \Ref{r=2}, since every term in
\Ref{HL} satisfies it. 

On the other hand, by Lemma \ref{dual lem}, the dimension of the space of
polynomials in $F_n^{(k,2)}$ of degree $d$ is equal to the dimension of the
subspace of degree $d$ in $E_n^{(k,2)}$. By Lemma \ref{span lem} it is bounded
from above by the number of admissible partitions in $\pi_n$ of weight $d$.

Since we have as many linearly independent function $P_\la$ of degree
$d$ as the number of admissible partitions in $\pi_n$ of weight $d$,
they form a basis and the bounds are actually equalities.
\end{proof}

\subsection{The character $\chi_{k,2}$}
>From Proposition \ref{r=2 prop} we have the following immediate corollary.

\begin{prop}\label{r=2 char}
The character of the space $F^{(k,2)}$ is given by
\bean\label{char 2}
\chi_{k,2}(z,v)=\prod_{s=0}^\infty
(1+v^sz+v^{2s}z^2+\dots+v^{ks}z^k)=\prod_{s=0}^\infty \frac{1-(v^sz)^{k+1}}{1-v^sz}.
\eean 
${}\qquad \qquad \qquad\qquad \qquad \qquad\qquad \qquad \qquad\qquad
\qquad \qquad \qquad \qquad \qquad\qquad \qquad $ $\square$
\end{prop}

Note that according to \cite{FJMM1} for generic $t$ a basis in the
space of functions satisfying the wheel condition related to the wheel
set $\{t,t^{-k-1}\}$ is given by Macdonald polynomials
parametrized by $(k,2,n)$-admissible partitions. The
corresponding character cannot be written in such a simple factored form.

Now from the comparison of dimensions we also obtain:

\begin{cor}\label{basis cor}
The elements $\{\bar e_\la\}$ with admissible partitions $\la\in\pi_n$
form a basis of $E_n^{(k,2)}$. $\square$
\end{cor}

Define formal power series in $v$, $b^{(k)}_n(v)$, by
\be
\chi_{k,2}(z,v)=\sum_{n=0}^\infty b_n^{(k)}(v)z^n.
\ee

For $k=1$, we have the following result.

\begin{lem}\label{1,2} We have 
\be
b^{(1)}_n(v)=\frac{v^{n(n-1)/2}}{\prod_{i=1}^n(1-v^i)}.
\ee
\end{lem}
\begin{proof}
The lemma follows from the identities
\be
\prod_{s=0}^\infty(1+v^sz)=\sum_{n=0}^\infty
\left(\sum_{\la}v^{|\la|}\right)z^n=
\sum_{n=0}^\infty \frac{v^{n(n-1)/2}}{\prod_{i=1}^n(1-v^i)}z^n,
\ee
where in the second expression
the sum is over all partitions $\la\in\pi_n$ such that
$\la_i>\la_{i+1}$, $i=1,\dots,n-1$. 
\end{proof}

Now we give a formula for the coefficients $b_n^{(k)}$.

\begin{lem} We have
\bean\label{coef}
b_n^{(k)}(v)=\sum_{a,b;\; (k+1)a+b=n}\frac{(-1)^a v^{(k+1)a(a-1)/2}}
{\prod_{i=1}^a(1-v^{(k+1)i})\prod_{j=1}^b(1-v^j)}.
\eean 
\end{lem}
\begin{proof}
We have the identity
\be
\frac1{\prod_{s=0}^\infty (1-v^sz)}=\sum_{b=0}^\infty \left(
\sum_{\la\in\pi_b}v^{|\la|}\right)z^b=\sum_{b=0}^\infty
\frac{z^b}{\prod_{j=1}^b(1-v^j)}.
\ee
We obtain the lemma multiplying this identity by the identity
\be
\prod_{s=0}^\infty(1-(v^sz)^{k+1})=\chi_{1,2}(-z^{k+1},v^{k+1})=
\sum_{a=0}^\infty
\frac{(-1)^a v^{(k+1)a(a-1)/2}z^{(k+1)a}}{\prod_{i=1}^a(1-v^{(k+1)i})},
\ee
which follows from Lemma \ref{1,2}.
\end{proof}

\subsection{Other bases in $F^{(k,2)}$}\label{mac basis sec}
The following more general proposition is proved by the methods 
of \cite{FJMM2}, Theorem 2.4. We skip the details of the proof and 
do not use this result in any other part of the paper.

\begin{prop}\label{r=2 general}
Let $t$ be a primitive root of unity of order $k+1$.
The set of Macdonald polynomials $\{P(x_1,\dots,x_n;t,q)\}$ where
$\la$ ranges over all admissible partitions is a basis of
$F_n^{(k,2)}$ if $q$ is not a root of unity. $\square$
\end{prop}

Proposition \ref{r=2 prop} is just $q=0$ case of Proposition \ref{r=2
  general}. 

\section{The case $r>2$}
\subsection{Frobenius homomorphism} Fix $k\in\Z_{>0}$, $r\in\Z_{>1}$ such that
  $k+1$ and $r-1$ are relatively prime. Fix 
primitive roots of unity $t$ and $q$ of order $k+1$ and $r-1$,respectively.
We consider the wheel set $\mc S_{r-1}(q,t)=\{t, tq,\dots, tq^{r-2}\}$.
It is invariant under the multiplication by $q$.

The space $F^{(k,r)}$ is the space of symmetric polynomials
$f(x_1,\dots,x_n)$ which vanish at
\bean\label{kr cond}
x_i=t^{i-1}q^{s_{i}}x_1 \qquad (1<i\leq k+1),
\eean
for all $s_2,\dots,s_{k+1}\in\{0,1,\dots,r-2\}$.
  
Define the {\it Frobenius homomorphism} 
\bea
\mc F:\;&& \C[x_1,\dots,x_n]\to \C[x_1,\dots,x_n],\\
&& f(x_1,\dots,x_n)\mapsto f(x_1^{r-1},\dots,x_n^{r-1}).
\eea

We have an obvious lemma.

\begin{lem} The Frobenius homomorphism induces an imbedding of the
  vector spaces
\be
\mc F:\;F^{(k,2)} \to F^{(k,r)}.
\ee
\qquad \qquad \qquad \qquad \qquad \qquad \qquad \qquad \qquad
\qquad\qquad \qquad \qquad \qquad 
\qquad \qquad $\square$
\end{lem}
We denote the image $\mc F(F^{(k,2)})$ by $R$.

Let $F_n^{(k,r)}\subset F^{(k,r)}$ be the subspace of polynomials with
 $n$ variables. Let $R_n=R\cap F_n^{(k,r)}$.

A symmetric polynomial $f\in F_n^{(k,r)}$ is in the Frobenius image 
$R_n\subset F_n^{(k,r)}$ 
if and only if there exists a symmetric 
polynomial $\tilde f(x_1,\dots,x_n)$ such that
$f(x_1,\dots,x_n)=\tilde f(x_1^{r-1},\dots,x_n^{r-1})$. 
In such a case $\tilde f\in F_n^{(k,2)}$. 

The space $F_n^{(k,r)}$ is an ideal in the algebra of symmetric polynomials.
We view the Frobenius image $R_n$ as a subring in
$\C[x_1,\dots,x_n]$. Our goal is to prove that $F_n^{(k,r)}$ is a free
$R$-module of rank $N=(r-1)^n$ 
and to describe the generating set (see Proposition \ref{separating}).

\subsection{Slim partitions} We describe the division with remainders
of partition by positive integers.

By definition, the result of addition of two partitions
$\la,\mu\in\pi_n$ 
is a partition $(\la+\mu)\in\pi_n$ with components
$\la_i+\mu_i$, $i=1,\dots,n$. In particular,
the result of multiplication of a partition $\la$ by a
positive integer $r-1$ is the partition $(r-1)\la$ with components
$(r-1)\la_i$.
A partition $\la$ is called divisible by a positive integer $r-1$ 
if all parts $\la_i$ are divisible by $r-1$. 


Let a symmetric polynomial $f(x_1,\dots,x_n)$ be
expressed in terms of monomial functions  
$f=\sum_\la c_\la m_\la$, $c_\la\in \C$. Then
$f(x_1,\dots,x_n)=\tilde f(x_1^{r-1},\dots,x_n^{r-1})$ for some symmetric
polynomial $\tilde f$ if and only if $c_\la =0 $ for all
$\la$ which are not divisible by $r-1$.

A partition $\la\in\pi_n$ is called {\it $(r-1)$-slim} if
$\la_i-\la_{i+1}<r-1$ for $i=1,\dots,n-1$ and $\la_n<r-1$. 
We denote the subset of $\pi_n$ consisting of all
$(r-1)$-slim partitions by $\pi_n^{s}$.
The cardinality of $\pi_n^{s}$ is $N=(r-1)^n$.

The following lemma is straightforward.

\begin{lem}
Let $\la\in\pi_n$. There is a unique way to represent $\la$ in the form
\be
\la=(r-1)\mu +\nu,
\ee
where $\mu,\nu\in\pi_n$ are partitions and $\nu$ is $(r-1)$-slim. $\square$
\end{lem}

The partitions $\mu, \nu$ in the lemma are called the {\it quotient}
and the {\it remainder} of partition $\la$ divided by integer $r-1$.

\subsection{Wheel condition is broken if the highest partition is slim}
A symmetric polynomial $f(x_1,\dots,x_n)$ is called of highest
partition $\la$ if it has the form
\be
f=c_\la m_\la+\sum_{\mu<\la}c_\mu m_\mu \qquad (c_\mu\in \C),
\ee
where $c_\la$ is a nonzero complex number.

Let $\mu\in\pi_n$ be an $(r-1)$-slim partition and let $g_\mu$ be of
highest partition $\mu$.
We claim that $g$ does not satisfy the wheel condition, 
that is $g\not\in F^{(k,r)}$.
We prove a slightly stronger statement.
\begin{prop}\label{broken}
Let $\mu\in\pi_n^s$ be an $(r-1)$-slim partition and let $g_\mu$ be of
highest partition $\mu$. If $n>k$ then for any non-zero complex number
$c$ there exist 
$s_1,\dots,s_k\in\{0,1,\dots,r-2\}$ such that
$g(cq^{s_1},ctq^{s_2},\dots,ct^kq^{s_{k+1}},x_{k+2},\dots,x_n)\neq 0$.
\end{prop}
\begin{proof}
We fix a non-zero complex number $c$.
Let $\C^N$ ($N=(r-1)^n$) be the vector space
with basis $\{u_s\}$ where the index $s=(s_1,\dots,s_n)$ ranges over all
sequences such that $s_i\in\{0,\dots,r-2\}$.
Let $m=n-k-1\geq 0$. For $y=(y_1,\dots,y_m)\in \C^m$, define the evaluation map
\bea
\kappa_y:\;&& \C[x_1,\dots,x_n]^{{\mathfrak S}_n} \to \C^N \\
&& f\mapsto
\sum_{s}f(cq^{s_1},ctq^{s_2},\dots,ct^kq^{s_{k+1}},y_1q^{s_{k+2}},
\dots,y_mq^{s_n})u_s.
\eea
We need to show that there exists $y$ such that $\kappa_y(g_\mu)\neq 0$. 

If $f\in F_n^{(k,r)}$ then $\kappa_y(f)=0$ for all $y$.
If $f$ is in the image of the Frobenius homomorphism then 
$\kappa_y(f)=a(y)\sum_s u_s$ for some $a(y)\in \C$.

For each $\la\in\pi_n$ choose a symmetric polynomial $f_\la$ with
highest partition $\la$. Then $\{f_\la\}_{\la\in\pi_n}$ is a basis in
$\C[x_1,\dots,x_n]^{{\mathfrak S}_n}$. For each slim partition
$\mu\in\pi_n^s$ choose a
symmetric polynomial $g_\mu$ with highest partition $\mu$. Then
$\{g_\mu \mc F(f_\la)\}_{\la\in\pi_n,\mu\in\pi_n^s}$ is also a basis
in  $\C[x_1,\dots,x_n]^{{\mathfrak S}_n}$. 

In particular the algebra
$\C[x_1,\dots,x_n]^{{\mathfrak S}_n}$ is a free module of rank $N$ with generators
$g_\mu$ over the image of the Frobenius homomorphism . Therefore,
if $\kappa_y(g_\mu)=0$ for some slim partition $\mu$ and all $y$, then for
all $y$ the map $\kappa_y$ is not surjective. We claim it is impossible. 

Indeed, for generic choice of $(y_1,\dots,y_m)\in\C^m$, 
the images of the $N$ evaluation points
$(cq^{s_1},ctq^{s_2},\dots,ct^kq^{s_{k+1}},y_1q^{s_{k+2}},\dots, y_mq^n)$
in the quotient space $A=\C^n/{\mathfrak S}_n$ are all distinct.
Therefore the algebra of polynomials
$\C[A]=\C[x_1,\dots,x_n]^{{\mathfrak S}_n}$
separates these points and the evaluation map is surjective.
\end{proof}

\subsection{The main results}
We still have $t$ and $q$ as in (\ref{roots}).
Recall Macdonald's operators $D_n^r$.
\begin{lem}\label{mac lem}
If $f\in F_n^{(k,r)}$ satisfies the wheel condition then $D_n^r(q,\tilde
t)f$ satisfies the wheel condition for all $r\in\Z_{>0}$ and $\tilde t\in \C$.
\end{lem}
\begin{proof}
If a function (not necessary symmetric) $f$ 
satisfies the wheel conditions then $T_{q,x_i}f$ also satisfies the
wheel conditions for all $i=1,\dots,n$. The statement follows from
this observation.
\end{proof}

Consider a symmetric polynomial $h(x_1,\dots,x_n)$ of the form
\bean\label{split form}
h(x_1,\dots,x_n)=\sum_{i=1}^s f_ig_i,
\eean
where the highest partitions of $g_i$ are all slim and distinct 
and $f_i=\mc F(\tilde f_i)$ are all in the image of the Frobenius
homomorphism.

\begin{prop}\label{separating}
The polynomial $h$ is in $F_n^{(k,r)}$ if and only if $f_i\in R_n$,
i.e., if and only if $\tilde f_i$ are in  $F_n^{(k,2)}$.
\end{prop}
\begin{proof}
The ``if'' part of the Proposition is obvious.

To prove the ``only if'' part, we
first assume that $g_i=P_{\la^{(i)}}(q,\tilde t)$ are Macdonald
polynomials with some $\tilde t\in \C$ which is not a root of unity and all
$\la^{(i)}$ are distinct and slim. Notice 
that the Macdonald polynomials are well defined.

We have
\bea
D_n(X;q,\tilde t)h= \sum_{i=1}^s D_n(X;q,\tilde t)(g_if_i)= \\
\sum_{i=1}^s
f_i D_n(X;q,\tilde t)(g_i)=
\sum_{i=1}^s\prod_{j=1}^n(1+Xq^{\la_j^{(i)}}\tilde t^{n-j}) f_ig_i.
\eea
also satisfies the wheel condition by Lemma \ref{mac lem}. Since 
eigenvalues  $\prod_{j=1}^n(1+Xq^{\la_j^{(i)}}\tilde t^{n-j})$,
$i=1,\dots,s$, are all distinct, we conclude that $g_if_i$ satisfies
the wheel condition for all $i$. Now, since for $g_i$ the wheel condition is
broken by Proposition \ref{broken}, we get $f_i$ is zero on one of
the vanishing planes and therefore $\tilde f_i\in F_n^{(k,2)}$. 


Now we prove the general case by induction on $s$.
Let $\la^{(i)}$ be the highest 
partitions of $g_i$. Without loss of generality we assume that
$\la^{(s)}$ is maximal in the set $\{\la^{(i)}\}$. Then $h$ can be
written in the form 
\be
h=a P_{\la^{(s)}}(q,\tilde t)f_s+\sum_\mu P_\mu(q,\tilde t) f_\mu
\qquad (a\in \C, \; a\neq 0),
\ee
where $f_\mu=\mc F(\tilde f_\mu)$ are in the image of Frobenius
homomorphism and the sum is over slim partitions $\mu$ such that
$\mu\neq \la^{(s)}$. 
Therefore, using the previous argument we conclude that
$f_s\in R_n$. In particular $\sum_{i=1}^{s-1}f_ig_i=h-f_sg_s\in F_n^{(k,r)}$ and by the
induction hypothesis we obtain $f_i\in R_n$ for all $i$.
\end{proof}

Let $g_\mu$, $\mu\in\pi_n^s$, be any polynomials with highest partitions
$\mu$. Let  $\tilde f_\la$, where $\la$ ranges over
admissible partitions, be a basis in $F_n^{(k,2)}$. We described some of
such bases in Sections \ref{HL basis sec}, \ref{mac basis sec}.
Let $f_\la={\mc F}(\tilde f_\la)$.

We have proved the following theorem.
\begin{theorem}\label{basis thm} 
The polynomials $\{f_\la g_\mu \}$, where $\mu$ ranges over all
$(r-1)$-slim partitions and $\la$ ranges over
admissible partitions, form a basis in $F_n^{(k,r)}$.$\square$.
\end{theorem}

As a corollary we obtain a description of the
character $\chi_{k,r}(z,v)$ of the space
$F^{(k,r)}$.

\begin{theorem}\label{char thm}
The character of the space $F^{(k,r)}$ is given by
\be
\chi_{k,r}(z,v)=\sum_{n=0}^\infty \left(b_n^{(k)}(v^{r-1})\prod_{s=1}^n\frac{1-v^{s(r-1)}}{1-v^s}\right)\;z^n,  
\ee
where the coefficients $b_n(v)$ are given by \Ref{coef}.
\end{theorem}
\begin{proof}
The $z^n$ term of the character $\chi_{k,r}(z,v)$  is 
the product of the character of $(r-1)$-slim partitions of length at most $n$
and the $z^n$ term of the character $\chi_{k,2}(z,v^{r-1})$. The character
$\chi_{k,2}(z,v)$ is computed in
Proposition \ref{r=2 char}, it's $z_n$ term is $b_n^{(k)}(v)$ and
the character of the $(r-1)$-slim partitions in $\pi_n$ is equal to the product
$\prod_{s=1}^n(1-v^{s(r-1)})/(1-v^s)$.
\end{proof}

\bigskip

\noindent
{\it Acknowledgments.}\quad This work is partially supported by the
Grant-in-Aid for Scientific Research (B2) no.12440039, no.14340040 and
(A1) no.13304010, Japan Society for the Promotion of Science. The work
of EM is partially supported by NSF grant DMS-0140460. 
YT is supported by the Japan Society for the Promotion of Science.

\end{document}